\documentclass[12pt]{article}
\usepackage{graphicx}
\usepackage{color}
\catcode`\@=11 \@addtoreset{equation}{section}

\catcode`\@=12

\newtheorem{Theorem}{Theorem}[section]
\newtheorem{Proposition}{Proposition}[section]
\newtheorem{Lemma}{Lemma}[section]
\newtheorem{Corollary}{Corollary}[section]
\newtheorem{Remark}{Remark}[section]

\newcommand{\bTheorem}[1]{
\begin{Theorem} \label{T#1} }
\newcommand{\eT}{\end{Theorem}}

\newcommand{\bProposition}[1]{
\begin{Proposition} \label{P#1}}
\newcommand{\eP}{\end{Proposition}}

\newcommand{\bLemma}[1]{
\begin{Lemma} \label{L#1} }
\newcommand{\eL}{\end{Lemma}}

\newcommand{\bCorollary}[1]{
\begin{Corollary} \label{C#1} }
\newcommand{\eC}{\end{Corollary}}

\newcommand{\bRemark}[1]{
\begin{Remark} \label{R#1} }
\newcommand{\eR}{\end{Remark}}

\newcommand{\bFormula}[1]{
\begin{equation} \label{#1}}
\newcommand{\eF}{\end{equation}}

\newcommand{\Ov}[1]{\overline{#1}}

\newcommand{\DC}{C^\infty_c}

\newcommand{\vr}{\varrho}

\newcommand{\vt}{\vartheta}
\newcommand{\vu}{\vc{u}}
\newcommand{\vc}[1]{{\bf #1}}
\newcommand{\qed}{\bigskip \rightline {Q.E.D.} \bigskip}
\newcommand{\Div}{{\rm div}_x}
\newcommand{\Grad}{\nabla_x}

\newcommand{\tn}[1]{\mbox {\F #1}}
\newcommand{\dx}{{\rm d} {x}}
\newcommand{\dt}{{\rm d} t }

\newcommand{\dxdt}{\dx \ \dt}
\newcommand{\intO}[1]{\int_{\Omega} #1 \ \dx}

\newcommand{\bProof}{{\bf Proof: }}

\newcommand{\ep}{\varepsilon}

\font\F=msbm10 scaled 1000

\definecolor{Cgrey}{rgb}{0.85,0.85,0.85}
\definecolor{Cblue}{rgb}{0.50,0.85,0.85}
\definecolor{Cred}{rgb}{1,0,0}
\definecolor{fancy}{rgb}{0.10,0.85,0.10}

\newcommand\Cbox[2]{%
    \newbox\contentbox%
    \newbox\bkgdbox%
    \setbox\contentbox\hbox to \hsize{%
        \vtop{
            \kern\columnsep
            \hbox to \hsize{%
                \kern\columnsep%
                \advance\hsize by -2\columnsep%
                \setlength{\textwidth}{\hsize}%
                \vbox{
                    \parskip=\baselineskip
                    \parindent=0bp
                    #2
                }%
                \kern\columnsep%
            }%
            \kern\columnsep%
        }%
    }%
    \setbox\bkgdbox\vbox{
        \color{#1}
        \hrule width  \wd\contentbox %
               height \ht\contentbox %
               depth  \dp\contentbox
        \color{black}
    }%
    \wd\bkgdbox=0bp%
    \vbox{\hbox to \hsize{\box\bkgdbox\box\contentbox}}%
    \vskip\baselineskip%
}

\date{}

\begin{document}


\title{On the weak solutions to the equations of a compressible heat conducting gas}

\author{Elisabetta Chiodaroli  \and Eduard Feireisl \thanks{The research of E.F. leading to these results has received funding from the European Research Council under the European Union's Seventh Framework
Programme (FP7/2007-2013)/ ERC Grant Agreement 320078.} \and Ond\v rej Kreml\thanks{The work of O.K. is part of the SCIEX project 11.152}}

\maketitle

\centerline{Mathematisches Institut, Universit\"at Leipzig}

\centerline{Augustusplatz 10, D-04009 Leipzig, Germany}

\bigskip

\centerline{Institute of Mathematics of the Academy of Sciences of the Czech Republic}

\centerline{\v Zitn\' a 25, 115 67 Praha 1, Czech Republic}

\bigskip

\centerline{Institut f\"ur Mathematik, Universit\"at Z\"urich}

\centerline{Winterthurerstrasse 190, CH-8057  Z\"urich, Switzerland}






\maketitle

\bigskip





\begin{abstract}

We consider the weak solutions to the Euler-Fourier system describing the motion of a compressible heat conducting gas. Employing the method
of convex integration, we show that the problem admits infinitely many global-in-time weak solutions for any choice of smooth initial data.
We also show that for any initial distribution of the density and temperature, there exists an initial velocity such that the associated initial-value problem possesses infinitely many solutions that conserve the total energy.

\end{abstract}


\section{Introduction}
\label{i}

The concept of \emph{weak solution} has been introduced in the mathematical theory of systems of (nonlinear) hyperbolic conservation laws to incorporate the inevitable singularities in their solutions that may develop in a finite time no matter how smooth and small the data are. As is well known, however, many nonlinear problems are not well posed in the weak framework and several classes of admissible weak solutions have been identified to handle
this issue. The implications of the \emph{Second law of thermodynamics} have been widely used in the form of various entropy conditions in order to
identify the physically relevant solutions. Although this approach has been partially successful when dealing with systems in the simplified
$1D-$geometry, see Bianchini and Bressan \cite{BiaBre}, Bressan \cite{BRESSAN}, Dafermos \cite{D4}, Liu \cite{LiuTP}, among others, the more realistic problems in higher spatial dimensions seem to
be out of reach of the theory mostly because the class of ``entropies'' is rather poor consisting typically of a single (physical) entropy. Recently,
De Lellis and Sz\' ekelyhidi \cite{DelSze} developed the method of \emph{convex integration} (cf. M\" uller and \v Sver\' ak \cite{MulSve}) to identify a large class of weak solutions to the Euler system
violating the principle of well-posedness in various directions. Besides the apparently non-physical solutions producing the kinetic energy
(cf. Shnirelman \cite{Shn}), a large
class of data has been identified admitting infinitely many weak solutions that comply with a major part of the known admissibility criteria, see
De Lellis and Sz\' ekelyhidi \cite{DelSze3}.

In this paper, we develop the technique of  \cite{DelSze3} to examine the well-posedness of the full
\emph{Euler-Fourier system}:

\Cbox{Cgrey}{

\bFormula{i1}
\partial_t \vr + \Div (\vr \vu) = 0,
\eF
\bFormula{i2}
\partial_t (\vr \vu) + \Div (\vr \vu \otimes \vu) + \Grad p(\vr,\vt) = 0,
\eF
\bFormula{i3}
\partial_t (\vr e(\vr, \vt) ) + \Div (\vr e(\vr,\vt) \vu ) + \Div \vc{q} = - p(\vr,\vt) \Div \vu,
\eF

}

\noindent where $\vr(t,x)$ is the mass density, $\vu = \vu(t,x)$ the velocity field, and $\vt$ the (absolute) temperature of a compressible, heat conducting gas, see Wilcox \cite{WIL}. For the sake of simplicity, we restrict ourselves to the case of perfect monoatomic gas, for which the pressure $p(\vr,\vt)$ and
the specific internal energy $e(\vr,\vt)$ are interrelated through the constitutive equations:
\bFormula{i3a}
p(\vr,\vt) = \frac{2}{3} \vr e(\vr,\vt), \ p(\vr,\vt) = a \vr \vt, \ a > 0.
\eF

Although the system (\ref{i1} - \ref{i3}) describes the motion in the absence of viscous forces, we suppose that the fluid is heat conductive, with
the heat flux $\vc{q}$ determined by the standard Fourier law:
\bFormula{i4}
\vc{q} = -\kappa \Grad \vt, \ \kappa > 0.
\eF

The problem  (\ref{i1} - \ref{i3}) is supplemented with the initial data
\bFormula{i5}
\vr(0, \cdot) = \vr_0, \ (\vr \vu)(0, \cdot) = \vr_0 \vu_0, \ \vt(0, \cdot) = \vt_0 \ \mbox{in} \ \Omega.
\eF
In addition, to avoid the effect of the kinematic boundary,
we consider the periodic boundary conditions, meaning the physical domain $\Omega$ will be taken the flat torus
\[
\Omega = \tn{T}^3 = \left( [ 0 , 1 ] |_{ \{ 0 , 1 \} } \right)^3.
\]

The first part of the paper exploits the constructive aspect of convex integration. We present a ``variable coefficients'' variant of a result
of De Lellis-Sz\' ekelyhidi \cite{DelSze3} and show the existence of \emph{infinitely many} global-in-time\footnote{By global-in-time solutions we mean here solutions defined on $[0,T)$ for any $T>0$. For discussion about solutions defined on $[0,\infty)$ see Section \ref{CCL}.} weak solutions to the problem
(\ref{i1} - \ref{i5}) for any physically relevant choice of (smooth) initial data. Here, physically relevant means that the initial distribution of
the density $\vr_0$ and the temperature $\vt_0$ are strictly positive in $\Omega$. These solutions satisfy also the associated entropy equation;
whence they comply with the Second law of thermodynamics.

Similarly to their counterparts constructed in \cite{DelSze3}, these ``wild'' weak solutions violate the First law of thermodynamics, specifically,
the total energy at any positive time is strictly larger than for the initial data. In order to eliminate the non-physical solutions, we therefore
impose the total energy conservation in the form:
\bFormula{i5a}
E(t) = \intO{ \vr \left( \frac{1}{2}  |\vu|^2 + e(\vr,\vt) \right)(t, \cdot) } =
\intO{ \vr_0 \left( \frac{1}{2}  |\vu_0|^2 + e(\vr_0,\vt_0) \right) } = E_0 \ \mbox{for (a.a.)}\ t \in (0,T).
\eF

Following \cite{FeiNov10} we show that the system (\ref{i1} - \ref{i3}), augmented with the total energy balance (\ref{i5a}), satisfies the principle of \emph{weak-strong uniqueness}. Specifically, the weak and strong solutions emanating from the same initial data necessarily coincide as long as the latter exists.
In other words, the strong solutions are unique in the class of weak solutions. This property remains valid even if we replace the internal energy
equation (\ref{i3}) by the entropy \emph{inequality}
\bFormula{i5b}
\partial_t (\vr s(\vr,\vt)) + \Div (\vr s(\vr,\vt)\vu ) + \Div \left( \frac{ \vc{q} }{\vt} \right) \geq - \frac{ \vc{q} \cdot \Grad \vt}{\vt^2},
\ \vt D s(\vr,\vt) \equiv D e(\vr,\vt) + p(\vr, \vt) D \left(\frac{1}{\vr} \right),
\eF
in the spirit of the theory developed in \cite{FeNo6}.

Although the stipulation of (\ref{i5a}) obviously eliminates the non-physical energy producing solutions, we will show that for
\emph{any} initial data $\vr_0$, $\vt_0$ there exists an initial velocity $\vu_0$ such that the problem (\ref{i1} - \ref{i5}) admits infinitely
many global-in-time weak solutions that satisfy the total energy balance (\ref{i5a}).

The paper is organized as follows. After a brief introduction of the concept of weak solutions in Section \ref{w}, we discuss the problem
of existence of infinitely many solutions for arbitrary initial data, see Section \ref{sec}. In particular, we prove a ``variable coefficients''
variant of a result of De Lellis and Sz\' ekelyhidi \cite{DelSze3} and employ the arguments based on Baire's category. In Section \ref{ec}, we show
the weak-strong uniqueness principle for the augmented system and then identify the initial data for which the associated
solutions conserve the total energy. The paper is concluded by some remarks on possible extensions in Section \ref{CCL}.

\section{Weak solutions}
\label{w}

To simplify presentation, we may assume, without loss of generality, that
\[
a = \kappa = 1.
\]

We say that a trio $[\vr,\vt,\vu]$ is a \emph{weak solution} of the problem (\ref{i1} - \ref{i5}) in the space-time cylinder $(0,T) \times \Omega$ if:

\begin{itemize}
\item
the density $\vr$ and the temperature $\vt$  are positive in $(0,T) \times \Omega$;
\item
\bFormula{w1}
\int_0^T \intO{ \left( \vr \partial_t \varphi + \vr \vu \cdot \Grad \varphi \right) } \ \dt = - \intO{ \vr_0 \varphi (0, \cdot) }
\eF
for any test function $\varphi \in \DC([0,T) \times \Omega)$;
\item
\bFormula{w2}
\int_0^T \intO{ \left( \vr \vu \cdot \partial_t \varphi + \vr \vu \otimes \vu : \Grad \varphi + \vr \vt \Div \varphi \right) } \ \dt
= - \intO{ \vr_0 \vu_0 \cdot \varphi(0, \cdot) }
\eF
for any test function $\varphi \in \DC([0,T) \times \Omega;R^3)$;
\item
\bFormula{w3}
\int_0^T \intO{ \left( \frac{3}{2} \left[ \vr \vt \partial_t \varphi + \vr \vt \vu \cdot \Grad \varphi \right] - \Grad \vt \cdot \Grad \varphi - \vr\vt\Div\vu\varphi \right) } \ \dt = - \intO{ \vr_0 \vt_0 \varphi (0, \cdot) }
\eF
for any test function $\varphi \in \DC([0,T) \times \Omega)$.
\end{itemize}

As a matter of fact, the weak solutions we construct in this paper will be rather regular with the only exception of the velocity field. In particular,
the functions $\vr$, $\vt$, and even $\Div \vu$ will be continuously differentiable in $[0,T] \times \Omega$, and, in addition,
\[
\vt \in L^2(0,T; W^{2,p}(\Omega)), \ \partial_t \vt \in L^p(0,T; L^p(\Omega)) \ \mbox{for any}\ 1 \leq p < \infty.
\]
Thus the equations (\ref{i1}), (\ref{i3}) will be in fact satisfied pointwise a.a. in $(0,T) \times \Omega$.
As for the velocity field, we have
\[
\vu \in C_{\rm weak}([0,T]; L^2(\Omega;R^3)) \cap L^\infty((0,T) \times \Omega;R^3),\ \Div \vu \in C([0,T] \times R^3).
\]

\section{Second law is not enough}
\label{sec}

Our first objective is to show the existence of infinitely many solutions to the Euler-Fourier system for arbitrary (smooth) initial data.

\Cbox{Cgrey}{

\bTheorem{sec1}
Let $T > 0$. Let the initial data satisfy
\bFormula{reg}
\vr_0 \in C^3(\Omega),
\vt_0 \in C^2(\Omega), \vu_0 \in C^3(\Omega;R^3), \ \vr_0(x) > \underline{\vr} > 0, \ \vt_0(x) > \underline{\vt} > 0 \ \mbox{for any}\
x \in \Omega.
\eF

Then the initial-value problem (\ref{i1} - \ref{i5}) admits infinitely many weak solutions in $(0,T) \times \Omega$ belonging to the class:
\[
\vr \in C^2([0,T] \times \Omega), \ \partial_t \vt \in L^p(0,T; L^p(\Omega)),\ \nabla^2_x \vt \in L^p(0,T; L^p(\Omega; R^{3 \times 3}))
\ \mbox{for any}\ 1 \leq p < \infty,
\]
\[
\vu \in C_{\rm weak}([0,T]; L^2(\Omega;R^3))
\cap L^\infty((0,T) \times \Omega ; R^3), \ \Div \vu \in C^2([0,T] \times \Omega).
\]

\eT

}

\bRemark{rem1}

Using the maximal regularity theory for parabolic equations (see Amann \cite{Amann1}, Krylov \cite{Krylov}) we observe that
$\vt$ is a continuous function of the time variable $t$ ranging in the interpolation space $[L^p(\Omega); W^{2,p}(\Omega)]_\alpha$ for any $1 \leq p < \infty$ finite and any
$\alpha \in (0,1)$. Thus it is possible to show that the conclusion of Theorem \ref{Tsec1} remains valid if we assume that
\[
\vt_0 \in [L^p(\Omega); W^{2,p}(\Omega)]_\alpha \ \mbox{for sufficiently large}\ 1\leq p < \infty\  \mbox{and}\ 0 < \alpha < 1, \ \vt_0 > 0 \ \mbox{in}\ \Omega,
\]
where $[,]_\alpha$ denotes the real interpolation.
In particular, the solution $\vt(t, \cdot)$ will remain in the same regularity class for any $t \in [0,T]$.

\eR

The rest of this section is devoted to the proof of Theorem \ref{Tsec1}.

\subsection{Reformulation}
\label{refor}

Following Chiodaroli \cite{Chiod}, we reformulate the problem in the new variables $\vr$, $\vt$, and $\vc{w} = \vr \vu$ obtaining, formally
\bFormula{sec1}
\partial_t \vr + \Div \vc{w} = 0,
\eF
\bFormula{sec2}
\partial_t \vc{w} + \Div \left( \frac{\vc{w} \otimes \vc{w} }{\vr} \right) + \Grad (\vr \vt)  = 0,
\eF
\bFormula{sec3}
\frac{3}{2} \left( \vr \partial_t \vt + \vc{w} \cdot \Grad \vt \right) = \Delta \vt - \vt \Div \vc{w} + \vt \frac{\Grad \vr}{\vr} \cdot \vc{w}.
\eF

Next, we take the following ansatz for the density:
\[
\varrho(t,x) = \vr_0(x) - h(t) \Div (\vr_0(x) \vu_0(x) ) \equiv \tilde \vr (t,x) ,
\]
with
\bFormula{sec3aa}
h \in C^2[0,T], \ h(0) = 0 , \ h' (0) = 1, \ \vr_0(x) - h(t) \Div (\vr_0(x) \vu_0 (x) ) > \frac{\underline{\vr}}{2}
\ \mbox{for all}\ t \in [0,T], \ x \in \Omega.
\eF

Accordingly, we write $\vc{w}$ in the form of its Helmholtz decomposition
\[
\vc{w} = \vc{v} + \Grad \Psi , \ \Div \vc{v} = 0,\  \Delta \Psi = h' (t) \Div (\vr_0 \vu_0), \ \intO{ \Psi } = 0.
\]
Obviously, by virtue of the hypotheses (\ref{reg}) imposed on the initial data, we have
\[
\tilde \vr \in C^2([0,T] \times \Omega), \Grad \Psi \in C^2([0,T] \times \Omega; R^3), \ \tilde \vr (0, \cdot) = \vr_0, \vc{w}_0 = \vr_0 \vu_0.
\]
Moreover,
the equation of continuity (\ref{sec1}) is satisfied pointwise in $(0,T) \times \Omega$, while the remaining two equations (\ref{sec1}), (\ref{sec2}) read \bFormula{sec2a}
\partial_t \vc{v} + \Div \left( \frac{(\vc{v} + \Grad \Psi) \otimes (\vc{v} + \Grad \Psi) }{\tilde \vr} \right) + \Grad \left( \tilde \vr \vt
+ \partial_t \Psi \right)  = 0,\ \Div \vc{v} = 0,
\eF
\bFormula{sec2aa}
\vc{v}(0, \cdot) = \vc{v}_0 = \vr_0 \vu_0 - \Grad \Delta^{-1}\Div ( \vr_0 \vu_0),
\eF
\bFormula{sec3a}
\frac{3}{2} \Big( \tilde \vr \partial_t \vt + (\vc{v} + \Grad \Psi)  \cdot \Grad \vt \Big) = \Delta \vt - \vt \Delta \Psi + \vt \frac{\Grad \tilde \vr}{\tilde \vr} \cdot (\vc{v} + \Grad \Psi), \ \vt(0, \cdot) = \vt_0.
\eF

\subsection{Internal energy and entropy equations}

For a given vector field $\vc{v} \in L^\infty((0,T) \times \Omega;R^3)$, the internal energy equation (\ref{sec3a}) is \emph{linear} with
respect to $\vt$ and as such admits a unique solution $\vt = \vt[\vc{v}]$ satisfying the initial condition
$\vt(0,\cdot) = \vt_0$. Moreover, the standard $L^p-$theory for parabolic equations (see e.g. Krylov \cite{Krylov}) yields
\[
\vt(t,x) > 0 \ \mbox{for all}\ t \in [0,T], \ x \in \Omega,
\]
\bFormula{i10-}
\partial_t \vt \in L^p(0,T; L^p(\Omega)), \ \nabla^2_x \vt \in L^p(0,T; L^p(\Omega; R^{3 \times 3}) \ \mbox{for any}\ 1 \leq p < \infty,
\eF
where the bounds depend only on the data and $\| \vc{v} \|_{L^\infty((0,T) \times \Omega;R^3)}$.

Dividing (\ref{sec3a}) by $\vt$ we deduce the \emph{entropy equation}
\bFormula{i10}
\tilde \vr \partial_t \log \left( \frac{\vt^{3/2}}{\tilde \vr} \right) + (\vc{v} + \Grad \Psi) \cdot \Grad \log \left( \frac{\vt^{3/2}}{\tilde \vr} \right)
= \Delta \log(\vt) + |\Grad \log(\vt) |^2,
\eF
where we have used the identity $- \Delta \Psi = \partial_t \tilde \vr$. We note that, given the regularity of the solutions
in Theorem \ref{Tsec1},
the entropy equation (\ref{i10}) and the internal energy equation (\ref{sec3a}) are \emph{equivalent}. In particular, the weak solutions we construct
are compatible with the Second law of thermodynamics.

\subsubsection{Uniform bounds}

Introducing a new variable
\[
Z = \log \left( \frac{\vt^{3/2}}{\tilde \vr} \right)
\]
we may rewrite (\ref{i10}) as
\bFormula{i11}
\tilde \vr \partial_t Z + \Big( \vc{v} + \Grad \Psi - \frac{8}{9} \Grad\log(\tilde \vr) \Big) \cdot \Grad Z = \frac{2}{3} \Delta Z + \frac{4}{9} |\Grad Z|^2 + \frac{2}{3} \Delta \log(\tilde \vr) +
\frac{4}{9} |\Grad \log(\tilde \vr) |^2.
\eF
Applying the standard parabolic comparison principle to (\ref{i11}) we conclude that $|Z|$ is bounded only in terms of the initial data and the time $T$. Consequently, the constants $\underline{\vt}$, $\Ov{\vt}$ can be taken in such a way that
\bFormula{i12}
0 < \underline{\vt} \leq \vt[\vc{v}] (t,x) \leq \Ov{\vt} \ \mbox{for all} \ t \in [0,T], \ x \in \Omega.
\eF
We emphasize that the constants $\underline{\vt}$, $\Ov{\vt}$ are \emph{independent} of $\vc{v}$ - a crucial fact that will be used in the future
analysis.

\subsection{Reduction to a modified Euler system}
\label{reduc}

Summing up the previous discussion, our task reduces to finding (infinitely many) solutions to the problem
\bFormula{sec2ab}
\partial_t \vc{v} + \Div \left( \frac{(\vc{v} + \Grad \Psi) \otimes (\vc{v} + \Grad \Psi) }{\tilde \vr} \right) + \Grad \left(  \tilde \vr \vt [ \vc{v}]
+ \partial_t \Psi - \frac{2}{3} \chi \right)  = 0,\ \Div \vc{v} = 0, \vc{v}(0, \cdot) = \vc{v}_0,
\eF
with a suitable spatially homogeneous function $\chi = \chi(t)$.

Following the strategy (and notation) of De Lellis and Sz\' ekelyhidi \cite{DelSze3}, we introduce the linear system
\bFormula{i14}
\partial_t \vc{v} + \Div \tn{U} = 0, \ \Div \vc{v} = 0,\  \vc{v}(0,\cdot) = \vc{v}_0, \ \vc{v}(T,\cdot) = \vc{v}_T ,
\eF
together with the function $\Ov{e}$,
\bFormula{i13}
\Ov{e}[\vc{v}] = \chi - \frac{3}{2} \tilde \vr \vt[\vc{v}] - \frac{3}{2} \partial_t \Psi,
\eF
with a positive function $\chi \in C[0,T]$ determined below.

Furthermore, we introduce the space $R^{3 \times 3}_{{\rm sym},0}$ of symmetric traceless matrices, with the operator norm
\[
\lambda_{\rm max} [ \tn{U} ] - \ \mbox{the maximal eigenvalue of}\ \tn{U} \in  R^{3 \times 3}_{{\rm sym},0}.
\]

Finally, we define the set of subsolutions
\bFormula{i15}
X_0 = \Bigg\{ \vc{v} \ \Big| \ \vc{v} \in L^\infty((0,T) \times \Omega;R^3) \cap C^1((0,T) \times \Omega;R^3)
\cap C_{\rm weak}([0,T]; L^2(\Omega;R^3)), 
\eF
\[
\vc{v} \ \mbox{satisfies (\ref{i14}) with some}\ \tn{U} \in C^1((0,T) \times \Omega ; R^{3 \times 3}_{{\rm sym},0}),
\]
\[
\inf_{t \in (\ep,T), x \in \Omega} \left\{ \Ov{e}[\vc{v}] - \frac{3}{2} \lambda_{\rm max} \left[ \frac{( \vc{v} + \Grad \Psi) \otimes (\vc{v} + \Grad \Psi) }{\tilde \vr} - \tn{U} \right] \right\} > 0 \ \mbox{for any} \ 0 < \ep < T \Bigg\}.
\]
Note that $X_0$ is different from its analogue introduced by Chiodaroli \cite{Chiod} and De Lellis and Sz\' ekelyhidi \cite{DelSze3}, in particular, the function
$\Ov{e}[\vc{v}]$ depends on the field $\vc{v}$.

As shown by De Lellis and Sz\' ekelyhidi \cite{DelSze3}, we have the (pointwise) inequality
\[
\frac{1}{2}  |\vc{w}|^2 \leq \frac{3}{2} \lambda_{\rm max} \left[ \vc{w} \otimes \vc{w} - \tn{U} \right], \ \vc{w} \in R^3,\  \tn{U} \in R^{3 \times 3}_{{\rm sym},0},
\]
where the identity holds only if
\[
\tn{U} = \vc{w} \otimes \vc{w} - \frac{1}{3} |\vc{w}|^2 \tn{I}.
\]
Consequently, by virtue of (\ref{i12}), there exists a constant $c$ depending only on the initial data $[\vr_0,\vt_0, \vu_0]$ such that
\bFormula{i16}
\sup_{t \in [0,T]} \| \vc{v}(t, \cdot) \|_{L^\infty(\Omega;R^3)} < c \ \mbox{for all} \ \vc{v} \in X_0.
\eF

Next, we choose the function $\chi \in C[0,T]$ in (\ref{i13}) so large that
\[
\frac{3}{2} \lambda_{\rm max} \left[ \frac{( \vc{v}_0 + \Grad \Psi) \otimes (\vc{v}_0 + \Grad \Psi) }{\tilde \vr}  \right]  < \chi - \frac{3}{2} \tilde \vr \vt[\vc{v}_0] - \frac{3}{2} \partial_t \Psi \equiv \Ov{e}[\vc{v}_0]
\ \mbox{for all} \ (t,x) \in [0,T] \times \Omega,
\]
in particular, the function $\vc{v}_0 = \vc{v}_0(x)$, together with the associated tensor $\tn{U} \equiv 0$, belongs to the set $X_0$, where $\vc{v}_0 = \vc{v}_T$.

We define a topological space $X$ as a completion of $X_0$ in $C_{\rm weak}([0,T]; L^2(\Omega;R^3))$ with respect to the metric $d$ induced by the weak topology of
the Hilbert space $L^2(\Omega;R^3)$.
As we have just observed, the space $X_0$ is non-empty as $\vc{v} = \vc{v}_0$ is in $X_0$.

Finally, we consider a family of functionals
\bFormula{i17}
I_\ep  [\vc{v}] = \int_\ep^T \intO{ \left( \frac{1}{2} \frac{|\vc{v} + \Grad \Psi|^2}{\tilde \vr} - \Ov{e} [\vc{v}] \right) } \ \dt \ \mbox{for}\ \vc{v} \in X,\ 0 < \ep < T.
\eF
As a direct consequence of the parabolic regularity estimates (\ref{i10-}), we observe that
\bFormula{i17a}
\Ov{e} [\vc{v}] \to \Ov{e} [\vc{w}] \ \mbox{in} \ C([0,T] \times \Omega) \ \mbox{whenever}\ \vc{v} \to \vc{w}\ \mbox{in}\  X;
\eF
therefore each $I_\ep$ is a compact perturbation of a convex functional; whence lower semi-continuous in $X$.

In order to proceed, we need the following crucial result that may be viewed as a ``variable coefficients'' counterpart of
\cite[Proposition 3]{DelSze3}.

\bProposition{w1}
Let $\vc{v} \in X_{0}$ such that
\[
I_\ep[\vc{v}] < - \alpha < 0,\ 0 < \ep < T/2.
\]

There there is $\beta = \beta(\alpha) > 0$ and a sequence $\{ \vc{v}_n \}_{n = 1}^\infty \subset X_0$ such that
\[
\vc{v}_n \to \vc{v} \ \mbox{in} \ C_{\rm weak}([0,T]; L^2(\Omega;R^3)),\
\liminf_{n \to \infty} I_\ep [\vc{v}_n]  \geq I_\ep [\vc{v}] + \beta.
\]
\eP

\noindent
We point out that the quantity $\beta = \beta(\alpha)$ is independent of $\ep$ and $\vc{v}$.

Postponing the proof of Proposition \ref{Pw1} to the next section, we complete the proof of Theorem \ref{Tsec1} following the line
of arguments of \cite{DelSze3}. To begin, we observe that cardinality of the space $X_0$ is infinite. Secondly, since each $I_\ep$ is a bounded lower semi-continuous functional on a complete metric space, the points of continuity of $I_\ep$ form a residual set in $X$. The set
\[
\mathcal{C} = \bigcap_{m > 1} \left\{ \vc{v} \in X \ | \ I_{1/m}[\vc{v}] \ \mbox{is continuous} \ \right\},
\]
being an intersection of a countable family of residual sets, is residual, in particular of infinite cardinality, see
De Lellis and Sz\' ekelyhidi \cite{DelSze3} for a more detailed explanation of these arguments.

Finally, we claim that for each $\vc{v} \in \mathcal{C}$ we have
\[
I_{1/m}[\vc{v}] = 0 \ \mbox{for all}\ m > 1;
\]
whence
\[
\frac{1}{2} \frac{|\vc{v} + \Grad \Psi|^2}{\tilde \vr} = \Ov{e}[\vc{v}] \equiv \chi - \frac{3}{2} \tilde \vr \vt[\vc{v}] - \frac{3}{2} \partial_t \Psi,
\]
\[
\tn{U} = \frac{ (\vc{v} + \Grad \Psi) \otimes (\vc{v} + \Grad \Psi) }{\tilde \vr} - \frac{1}{3} \frac{|\vc{v} + \Grad \Psi|^2}{\tilde \vr} \tn{I}
\ \mbox{for a.a.}\ (t,x) \in (0,T) \times \Omega,
\]
in other words, the function $\vc{v}$ is a weak solution to the problem (\ref{sec2ab}). Indeed, assuming $I_{1/m}[\vc{v}] < -2\alpha < 0$, we first find a sequence $\{ \vc{u}_n \}_{n=1}^\infty \subset X_0$ such that
\[
\vc{u}_n \to \vc{v} \ \mbox{in} \ C_{\rm weak}([0,T]; L^2(\Omega; R^3)), \ I_{1/m}[\vc{u}_n] < -\alpha.
\]
Then for each $\vc{u}_n$ we use
Proposition \ref{Pw1} and together with standard diagonal argument we obtain a sequence $\{ \vc{v}_n \}_{n=1}^\infty \subset X_0$ such that
\[
\vc{v}_n \to \vc{v} \ \mbox{in} \ C_{\rm weak}([0,T]; L^2(\Omega; R^3)), \ \liminf_{n \to \infty} I_{1/m}[\vc{v}_n] \geq I_{1/m}[\vc{v}] + \beta,
\ \beta > 0,
\]
in contrast with the fact that $\vc{v}$ is a point of continuity of $I_{1/m}$.

\subsection{Proof of Proposition \ref{Pw1}}

The proof of Proposition \ref{Pw1} is based on a localization argument, where variable coefficients are replaced by constants. The
fundamental building block is the following result proved by De Lellis and Sz\' ekelyhidi \cite[Proposition 3]{DelSze3}, Chiodaroli
\cite[Section 6, formula (6.9)]{Chiod}:

\bLemma{P1}
Let $[T_1,T_2]$, $T_1 < T_2$, be a time interval and $B \subset R^3$ a domain. Let $\tilde r \in (0,\infty) $, $\tilde \vc {V} \in R^3 $ be \emph{constant} fields such that
\[
0 < \underline{r} < \tilde r < \Ov{r}, | \tilde \vc{V} | < \Ov{V}.
\]
Suppose that
\[
\vc{v} \in C_{\rm weak} ([T_1, T_2]; L^2(B, R^3)) \cap C^1 ((T_1,T_2) \times \Ov{B}; R^3)
\]
satisfies the linear system
\[
\partial_t \vc{v} + \Div \tn{U} = 0, \ \Div \vc{v} = 0 \ \mbox{in}\ (T_1,T_2) \times B
\]
with the associated field $\tn{U} \in C^1((T_1,T_2) \times \Ov{B}; R^{3 \times 3}_{{\rm sym},0})$ such that
\[
\frac{3}{2} \lambda_{\rm max} \left[ \frac{(\vc{v} + \tilde \vc{V}) \otimes (\vc{v} + \tilde \vc{V})}{\tilde r} - \tn{U} \right] <
e \ \mbox{in}\ (T_1,T_2) \times B
\]
for a certain function $e \in C([T_1;T_2] \times \Ov{B})$.

Then there exist sequences $\{ \vc{w}_n \}_{n=1}^\infty \subset \DC((T_1,T_2) \times B; R^3)$,
$\{ \tn{Y}_n \}_{n=1}^\infty \subset \DC((T_1, T_2) \times B; R^{3 \times 3}_{{\rm sym},0} )$ such that $\vc{v}_n = \vc{v} + \vc{w}_n$, $\tn{U}_n =
\tn{U} + \tn{Y}_n$ satisfy
\[
\partial_t \vc{v}_n + \Div \tn{U}_n = 0, \ \Div \vc{v}_n = 0 \ \mbox{in}\ (T_1,T_2) \times B
\]
\[
\frac{3}{2} \lambda_{\rm max} \left[ \frac{(\vc{v}_n + \tilde \vc{V}) \otimes (\vc{v}_n + \tilde \vc{V})}{\tilde r} - \tn{U}_n \right] <
e  \ \mbox{in}\ (T_1,T_2) \times B,
\]
\[
\vc{v}_n \to \vc{v} \in C_{\rm weak}([T_1,T_2]; L^2(B; R^3)),
\]
and
\bFormula{lala}
\liminf_{n \to \infty} \int_{T_1}^{T_2} \int_B \left| \vc{v}_n - \vc{v} \right|^2 \ \dxdt \geq
\Lambda \left( \underline{r}, \Ov{r}, \Ov{V}, \| e \|_{L^\infty((T_1,T_2) \times B)} \right) \int_{T_1}^{T_2} \int_B \left( e - \frac{1}{2} \frac{| \vc{v} + \tilde \vc{V} |^2}{\tilde r} \right)^2
\ \dxdt .
\eF

\eL

\bRemark{w1}
Note that $\tilde \vc{V}$ is constant in Lemma \ref{LP1}; whence
\[
\partial_t \vc{v} = \partial_t (\vc{v} + \tilde \vc{V}).
\]
\eR

\bRemark{w2}
It is important that the constant $\Lambda$ depends only on the quantities indicated explicitly in (\ref{lala}), in particular $\Lambda$ is
independent of $\vc{v}$, of the length of the time interval, and of the domain $B$.
\eR

\subsubsection{Localization principle}

The scale invariance encoded in (\ref{lala}) can be used for showing a ``variable coefficients'' variant of Lemma \ref{LP1}, specifically
when both $\tilde r$ and $\tilde \vc{V}$ are sufficiently smooth
functions of $t$ and $x$.

\bLemma{P2}
Let $\tilde \vr \in C^1([T_1,T_2] \times \Omega)$, ${\vc V} \in C^1([T_1,T_2] \times \Omega ; R^3 )$, $T_1 < T_2$ be functions satisfying
\[
0 < \underline{r} < \tilde \vr (t,x)  < \Ov{r}, | \vc{V}(t,x) | < \Ov{V} \ \mbox{for all}\ t,x.
\]
Suppose that
\[
\vc{v} \in C_{\rm weak} ([T_1, T_2]; L^2(\Omega, R^3)) \cap C^1 ((T_1,T_2) \times \Omega; R^3)
\]
solves the linear system
\[
\partial_t \vc{v} + \Div \tn{U} = 0, \ \Div \vc{v} = 0 \ \mbox{in}\ (T_1,T_2) \times \Omega
\]
with the associated field $\tn{U} \in C^1((T_1,T_2) \times B; R^{3 \times 3}_{{\rm sym},0})$ such that
\bFormula{P1}
\frac{3}{2} \lambda_{\rm max} \left[ \frac{(\vc{v} + \vc{V}) \otimes (\vc{v} + \vc{V})}{\tilde \vr} - \tn{U} \right] <
e - \delta \ \mbox{in} \ (T_1,T_2) \times \Omega
\eF
for some $e \in C([T_1;T_2] \times \Ov{B})$ and $\delta > 0$.

Then there exist sequences $\{ \vc{w}_n \}_{n=1}^\infty \subset \DC((T_1,T_2) \times \Omega; R^3)$,
$\{ \tn{Y}_n \}_{n=1}^\infty \subset \DC((T_1, T_2) \times \Omega; R^{3 \times 3}_{{\rm sym},0} )$ such that $\vc{v}_n = \vc{v} + \vc{w}_n$, $\tn{U}_n =
\tn{U} + \tn{Y}_n$ satisfy
\bFormula{P1A}
\partial_t \vc{v}_n + \Div \tn{U}_n = 0, \ \Div \vc{v}_n = 0 \ \mbox{in}\ (T_1,T_2) \times \Omega
\eF
\bFormula{P1B}
\frac{3}{2} \lambda_{\rm max} \left[ \frac{(\vc{v}_n +  \vc{V}) \otimes (\vc{v}_n +  \vc{V})}{\tilde \vr} - \tn{U}_n \right] <
e  \ \mbox{in}\  C((T_1;T_2) \times \Omega ),
\eF
\bFormula{P1C}
\vc{v}_n \to \vc{v} \in C_{\rm weak}([T_1,T_2]; L^2(\Omega; R^3)),
\eF
and
\bFormula{P1D}
\liminf_{n \to \infty} \int_{T_1}^{T_2} \int_\Omega \left| \vc{v}_n - \vc{v} \right|^2 \ \dxdt \geq
\Lambda \left( \underline{r}, \Ov{r}, \Ov{V}, \| e \|_{L^\infty((T_1,T_2) \times \Omega)} \right) \int_{T_1}^{T_2} \int_\Omega \left( e - \frac{1}{2} \frac{| \vc{v} +  \vc{V} |^2}{\tilde \vr} \right)^2
\ \dxdt .
\eF

\eL

\bRemark{w3}
The role of the positive parameter $\delta$ in (\ref{P1}) is only to say that the inequality (\ref{P1}) is strict, otherwise the conclusion of the
lemma is independent of the specific value of $\delta$.
\eR

\bRemark{w4}
In view of (\ref{P1C}), the convergence formula (\ref{P1D}) may be \emph{equivalently} replaced by
\bFormula{P1Da}
\liminf_{n \to \infty} \int_{T_1}^{T_2} \int_\Omega \frac{1}{2} \frac{  \left| \vc{v}_n + \vc{V} \right|^2 }{\tilde \vr}  \ \dxdt
\eF
\[
\geq \int_{T_1}^{T_2} \int_\Omega \frac{1}{2} \frac{  \left| \vc{v} + \vc{V} \right|^2 }{\tilde \vr}  \ \dxdt +
\Lambda \left( \underline{r}, \Ov{r}, \Ov{V}, \| e \|_{L^\infty((T_1,T_2) \times \Omega)} \right) \int_{T_1}^{T_2} \int_\Omega \left( e - \frac{1}{2} \frac{| \vc{v} +  \vc{V} |^2}{\tilde \vr} \right)^2
\ \dxdt .
\]
\eR

\bProof

We start with an easy observation that there exists $\ep = \ep\left( \delta, |e| \right)$ such that
\bFormula{P2}
\left\{
\begin{array}{c}
\frac{3}{2} \left| \lambda_{\rm max} \left[ \frac{(\vc{v} + \vc{V}) \otimes (\vc{v} + \vc{V})}{\tilde \vr} - \tn{U} \right] -
\lambda_{\rm max} \left[ \frac{(\vc{v} + \tilde \vc{V}) \otimes (\vc{v} + \tilde \vc{V})}{\tilde r} - \tn{U} \right] \right| < \frac{\delta}{4},
\\ \\
\left| \frac{1}{2} \frac{| \vc{v} + \vc{V} |^2}{\tilde \varrho } - \frac{1}{2} \frac{| \vc{v} + \tilde \vc{V} |^2}{\tilde r } \right| < \frac{\delta}{4}
\end{array}
\right\}
\eF
whenever
\[
\frac{3}{2} \lambda_{\rm max} \left[ \frac{(\vc{v} + \vc{V}) \otimes (\vc{v} + \vc{V})}{\tilde \vr} - \tn{U} \right] < e, \ |\tilde \vr - \tilde r | < \ep, \ | \vc{V} - \tilde \vc{V} | < \ep .
\]

For $\delta$ appearing in (\ref{P1}), we fix $\ep = \ep \left(\delta, \| e \|_{L^\infty(T_1,T_2) \times \Omega)} \right)$ as in (\ref{P2}) and find a (finite) decomposition of the set
$(T_1,T_2) \times \Omega$ such that
\[
[T_1, T_2] \times \Omega = \cup_{i = 1}^N \Ov{Q}_i ,\ Q_i = (T^i_1,T^i_2) \times {B}_i, \ Q_i \cap Q_j = \emptyset \ \mbox{for}\ i \ne j,
\]
\[
\sup_{Q_i}  \tilde \vr - \inf_{Q_i} \tilde \vr < \ep, \ \sup_{Q_i} \left| \vc{V} - \frac{1}{|Q_i|} \int_{Q_i} \vc{V} \ \dxdt \right| < \ep,
\]
where the number $N$ depends on $\ep$ and the Lipschitz constants of $\tilde \vr$, $\vc{V}$ in $[T_1,T_2] \times \Omega$.

Now, we apply Lemma \ref{LP1} on each set $Q_i$ with the choice of parameters
\[
\tilde r = \sup_{Q_i} \tilde \vr, \ \tilde \vc{V} = \frac{1}{|Q_i|} \int_{Q_i} \vc{V} \ \dxdt .
\]
In accordance with (\ref{P1}), (\ref{P2}), we have
\[
\frac{3}{2} \lambda_{\rm max} \left[ \frac{(\vc{v} + \tilde \vc{V}) \otimes (\vc{v} + \tilde \vc{V} )}{\tilde r} - \tn{U} \right] <
e  - \frac{\delta}{2} \ \mbox{in}\ Q_i.
\]
Under these circumstances,
Lemma \ref{LP1} yields a sequence of smooth functions $\vc{v}^i_n$, $\tn{U}^i_n$, with $\vc{v} - \vc{v}^i_n$, $\tn{U} - \tn{U}^i_n$ compactly supported in $Q_i$, such that
\[
\partial_t \vc{v}^i_n + \Div \tn{U}^i_n = 0, \ \Div \vc{v}^i_n = 0 \ \mbox{in}\ Q_i,
\]
\bFormula{P3A}
\frac{3}{2} \lambda_{\rm max} \left[ \frac{(\vc{v}^i_n + \tilde \vc{V}) \otimes (\vc{v}^i_n + \tilde \vc{V} )}{\tilde r} - \tn{U}^i_n \right] <
e  - \frac{\delta}{2},
\eF
\[
\vc{v}^i_n \to \vc{v} \ \mbox{in}\ C_{\rm weak} ([T^i_1, T^i_2], L^2(B_i)),
\]
and
\bFormula{P5}
\liminf_{n \to \infty} \int_{Q_i} \left| \vc{v}^i_n - \vc{v} \right|^2 \ \dxdt \geq
\Lambda \left( \underline{r}, \Ov{r}, \Ov{V}, \| e \|_{L^\infty((T_1,T_2) \times \Omega)} \right) \int_{Q_i} \left( e - \frac{1}{2} \frac{| \vc{v} + \tilde \vc{V} |^2}{\tilde r} - \frac{\delta}{2} \right)^2
\ \dxdt .
\eF

In view of (\ref{P2}), we replace $\tilde r$ by $\tilde \vr$ and $\tilde \vc{V}$ by $\vc{V}$ in
(\ref{P3A}) to obtain
\[
\frac{3}{2} \lambda_{\rm max} \left[ \frac{(\vc{v}^i_n +  \vc{V}) \otimes (\vc{v}^i_n +  \vc{V} )}{\tilde \vr} - \tn{U}_n \right] <
e  \ \mbox{in}\ \Ov{Q}_i.
\]
As $\vc{v}_n$, $\tn{U}_n$ are compactly supported perturbations of $\vc{v}$, $\tn{U}$ in $Q_i$, we may
define
\[
\vc{v}_n (t,x) = \vc{v}^i_n(t,x), \ \tn{U}_n = \tn{U}_n^i \ \mbox{for any}\ (t,x) \in \Ov{Q}_i, \ i=1,\dots, N.
\]

In accordance with the previous discussion, $\vc{v}_n$, $\tn{U}_n$ satisfy (\ref{P1A} - \ref{P1C}). In
order to see (\ref{P1D}), use (\ref{P2}) to observe that
\[
\left( e - \frac{1}{2} \frac{| \vc{v} + \tilde \vc{V} |^2}{\tilde r} - \frac{\delta}{2} \right) >
\left( e - \frac{1}{2} \frac{| \vc{v} + \vc{V} |^2}{\tilde \vr} - \frac{3 \delta}{4} \right) > 0 \ \mbox{in}\ Q_i;
\]
whence, making use of the hypothesis (\ref{P1}), specifically of the fact that
\[
e - \frac{1}{2} \frac{| \vc{v} + \vc{V} |^2}{\tilde \vr} > \delta,
\]
we may infer that
\[
\left( e - \frac{1}{2} \frac{| \vc{v} + \vc{V} |^2}{\tilde \vr} - \frac{3 \delta}{4} \right) \geq
\frac{1}{4} \left( e - \frac{1}{2} \frac{| \vc{v} + \vc{V} |^2}{\tilde \vr} \right) \ \mbox{in}\ Q_i.
\]
Thus, summing up the integrals in (\ref{P5}) we get (\ref{P1D}).

\qed

\subsubsection{Application to functionals $I_\ep$}

Fixing $\ep \in (0, T/2)$ we complete the proof of Proposition \ref{Pw1}.
Given $e \in C([0,T] \times \Omega)$, we introduce the spaces
\bFormula{i18}
X_{0, e} = \Bigg\{ \vc{v} \ \Big| \ \vc{v} \in  C^1((0,T) \times \Omega;R^3)
\cap C_{\rm weak}([0,T]; L^2(\Omega;R^3))
\eF
\[
 \vc{v} \ \mbox{satisfies (\ref{i14}) with some}\ \tn{U} \in C^1((0,T) \times \Omega ; R^{3 \times 3}_{{\rm sym},0}),
\]
\[
\left.
\frac{3}{2} \lambda_{\rm max} \left[ \frac{( \vc{v} + \Grad \Psi) \otimes (\vc{v} + \Grad \Psi) }{\tilde \vr} - \tn{U} \right] < e
\ \mbox{for}\ t \in (0,T),\ x \in \Omega \right\},
\]
along with the associated functionals
\bFormula{i19}
I_{\ep, e} [\vc{v}] = \int_\ep^T \intO{ \left( \frac{1}{2} \frac{ |\vc{v} + \Grad \Psi|^2 }{\tilde \vr} - e \right) } \ \dt \ \mbox{for}\ \vc{v} \in X,\ 0 < \ep < T/2.
\eF

The following assertion is a direct consequence of Lemma \ref{LP2}.

\bLemma{i1}
Let $\vc{v} \in X_{0, e}$, $ e \in C([0,T]\times \Omega)$, $0 < \ep < T/2$ be such that
\[
I_{\ep, e} [\vc{v}] < - \alpha < 0.
\]

There there is $\beta = \beta(\alpha, \| e \|_{L^\infty((0,T) \times \Omega)} ) > 0$, independent of $\ep$, and a sequence $\{ \vc{v}_n \}_{n > 0} \subset X_{0, e}$ such that
\[
\vc{v}_n \equiv \vc{v} \ \mbox{in}\ [0, \ep] \times \Omega,
\]
\[
\vc{v}_n \to \vc{v} \ \mbox{in} \ C_{\rm weak}([0,T]; L^2(\Omega;R^3)),\
\liminf_{n \to \infty} I_{\ep, e}[\vc{v}_n]  \geq I_{\ep, e}[\vc{v}] + \beta.
\]
\eL

\bRemark{w5}
We have used Lemma \ref{LP2} with (\ref{P1Da}), where, by virtue of Jensen's inequality,
\[
\int_{\ep}^{T} \int_\Omega \left( e - \frac{1}{2} \frac{| \vc{v} +  \vc{V} |^2}{\tilde \vr} \right)^2
\ \dxdt \geq \frac{1}{(T- \ep)|\Omega| } \left( \int_{\ep}^{T} \int_\Omega \left( e - \frac{1}{2} \frac{| \vc{v} +  \vc{V} |^2}{\tilde \vr} \right)
\ \dxdt \right)^2 \geq  \frac{\alpha^2}{(T- \ep)|\Omega| }.
\]
\eR

Finally, we show how Lemma \ref{Li1} implies Proposition \ref{Pw1}. Under the hypotheses of Proposition \ref{Pw1} and
in accordance with the definition of the space $X_0$, we find $\delta > 0$ and a function $e \in C([0,T] \times \Omega)$ such that
\[
e \leq \Ov{e}[\vc{v}] , \ e \equiv \Ov{e}[\vc{v}] - \delta \ \mbox{whenever}\ t \in [\ep, T],
\]
and
\[
\vc{v} \in X_{0,e}.
\]
Thus, we have
\[
I_{\ep,e}[\vc{v}] = \int_\ep^T \intO{ \left( \frac{1}{2} \frac{ |\vc{v} + \Grad \Psi|^2 }{\tilde \vr} - \Ov{e}[\vc{v}] + \delta  \right) } \ \dt =
I_\ep [\vc{v}] + (T- \ep)|\Omega| \delta < - \alpha / 2 < 0
\]
as soon as $\delta > 0$ was chosen small enough.

Consequently, by virtue of Lemma \ref{Li1}, there is a sequence of functions $\{ \vc{v}_n \}_{n =1}^\infty$ and $\beta = \beta(\alpha) > 0$ such that
\[
\vc{v}_n \in X_{0,e}, \ \vc{v}_n \equiv \vc{v} \ \mbox{in}\ [0,\ep] \times \Omega,
\]
and
\[
\vc{v}_n \to \vc{v} \ \mbox{in}\ C_{\rm weak}([0,T]; L^2(\Omega;R^3)), \ \liminf_{n \to \infty} I_{\ep,e}[\vc{v}_n] \geq  I_{\ep, e}[\vc{v}] + \beta = I_{\ep}[\vc{v}]
+ \beta + (T- \ep)|\Omega| \delta.
\]
Moreover, in accordance with (\ref{i17a}),
\[
I_{\ep,e}[\vc{v}_n] - I_{\ep}[\vc{v}_n] = \int_\ep^T \intO{ \Ov{e}[\vc{v}_n] - \Ov{e}[\vc{v}] + \delta } \to (T- \ep)|\Omega| \delta
\ \mbox{as}\ n \to \infty;
\]
whence we may infer that
\[
\liminf_{n \to \infty} I_{\ep}[\vc{v}_n] \geq  I_{\ep}[\vc{v}] + \beta.
\]

Finally, it remains to observe that $\vc{v}_n \in X_0$ for all $n$ large enough. To this end, note that
\[
\frac{3}{2} \lambda_{\rm max} \left[ \frac{( \vc{v}_n + \Grad \Psi) \otimes (\vc{v}_n + \Grad \Psi)}{\tilde \vr} - \tn{U} \right] =
\frac{3}{2} \lambda_{\rm max} \left[ \frac{( \vc{v} + \Grad \Psi)  \otimes (\vc{v} + \Grad \Psi)}{\tilde \vr} - \tn{U} \right] < e \leq \Ov{e}[\vc{v}] = \Ov{e}[\vc{v}_n]
\]
for all $t \in [0, \ep]$,
while
\[
\frac{3}{2} \lambda_{\rm max} \left[ \frac{(\vc{v}_n + \Grad \Psi) \otimes (\vc{v}_n + \Grad \Psi) }{\tilde \vr} - \tn{U} \right] < e = \Ov{e}[\vc{v}] - \delta \leq
\Ov{e}[\vc{v}_n] - \delta/2 \ \mbox{for all}\ t \in [\ep, T]
\]
for all $n$ large enough. We have proved Proposition \ref{Pw1}.

\section{Dissipative solutions}
\label{ec}

The solutions of the Euler-Fourier system constructed in Section \ref{sec} suffer an essential deficiency, namely they do not comply with the
First law of thermodynamics, meaning, they violate the total energy conservation (\ref{i5a}). On the other hand, the initial data in (\ref{reg}) are smooth enough for the problem to possess a standard classical solution existing on a possibly short time interval $(0, T_{\rm max})$, see e.g. Alazard \cite{AL1}, \cite{AL}. Note that the Euler-Fourier system fits also in the general framework and the corresponding existence theory developed by
Serre \cite{Serr3}, \cite{Serr4}. As the classical solutions are unique (in their regularity class) and obviously satisfy the total energy balance (\ref{i5a}),
the latter can be added to (\ref{w1} - \ref{w3}) as an admissibility condition.  The weak solutions of (\ref{i1} - \ref{i5}) satisfying (\ref{i5a}) will be called \emph{dissipative solutions}.

\subsection{Relative entropy (energy) and weak-strong uniqueness}

Following \cite{FeiNov10} we introduce the relative entropy functional
\bFormula{ec1}
\mathcal{E} \left( \vr , \vt , \vu \Big| r , \Theta, \vc{U} \right)
\eF
\[
= \intO{ \left( \frac{1}{2} \vr |\vu - \vc{U} |^2 + {H_{\Theta}
(\vr, \vt)} - \frac{ \partial H_{\Theta} (r, \Theta)}{\partial \vr}
(\vr - r) - H_{\Theta} ( r , \Theta ) \right) },
\]
where $H_\Theta$ is the ballistic free energy,
\[
H_\Theta(\vr, \vt) = \vr \Big( e(\vr, \vt) - \Theta s(\vr, \vt) \Big) = \vr \left( \frac{3}{2} \vt - \Theta \log \left( \frac{\vt^{3/2}}{\vr} \right)
\right).
\]
Repeating step by step the arguments of \cite{FeiNov10} we can show that any dissipative solution of the problem (\ref{i1} - \ref{i5}) satisfies
the \emph{relative entropy inequality}:
\bFormula{ec2}
\left[ \mathcal{E} \left( \vr, \vt, \vu \Big| r, \Theta, \vc{U} \right) \right]_{t=0}^{t = \tau}
+ \int_0^\tau \intO{ \Theta
 \frac{ |\Grad \vt|^2 }{\vt^2} )} \ \dt
\eF
\[
\leq \int_0^\tau \intO{ \Big( \vr ( \vc{U} - \vu) \cdot \partial_t \vc{U} + \vr
(\vc{U} - \vu) \otimes \vu : \Grad \vc{U} - p(\vr, \vt) \Div \vc{U} \Big) } \ \dt
\]
\[
- \int_0^\tau \intO{ \left( \vr \Big( s(\vr, \vt) - s(r, \Theta) \Big) \partial_t \Theta + \vr
\Big( s(\vr, \vt) - s(r, \Theta) \Big) \vu \cdot \Grad \Theta \right) } \ \dt
\]
\[
+ \int_0^\tau \intO{ \left( \left( 1 - \frac{\vr}{r} \right) \partial_t p(r, \Theta) - \frac{\vr}{r} \vu \cdot
\Grad p(r, \Theta)  \right) } + \int_0^\tau \intO{ \frac{\Grad \vt}{\vt} \cdot \Grad
\Theta  }
\ \dt
\]
for any trio of smooth ``test'' functions
\[
r,\ \Theta, \ \vc{U},\ r > 0,\ \Theta > 0.
\]

We report the following result \cite[Theorem 6.1]{EF101}.

\Cbox{Cgrey}{

\bTheorem{ec1}{\bf [Weak-strong uniqueness]}

Let $[\vr, \vt, \vu]$ be a dissipative (weak) solution of the problem (\ref{i1} - \ref{i5}), emanating from the
initial data $[\vr_0, \vt_0, \vu_0]$ satisfying (\ref{reg}), such that
\[
0 < \underline{\vr} < \vr(t,x) < \Ov{\vr},\ 0 < \underline{\vt} < \vt(t,x) < \Ov{\vt},\
| \vu (t,x) | < \Ov{u} \ \mbox{for a.a.}\ (t,x) \in (0,T) \times \Omega.
\]
Suppose that the same problem (with the same initial data) admits a classical solution $[\tilde \vr, \tilde \vt, \tilde \vu]$ in
$(0,T) \times \Omega$.

Then
\[
\vr \equiv \tilde \vr,\ \vt \equiv \tilde \vt,\ \vu \equiv \tilde \vu.
\]

\eT

}

\bRemark{ec1}
Here, ``classical'' means that all the necessary derivatives appearing in the equations are continuous functions in $[0,T] \times \Omega$.
\eR

\bRemark{ec2}
The proof of Theorem \ref{Tec1} is based on taking $r = \tilde \vr$, $\Theta = \tilde \vt$, $\vc{U} = \tilde \vu$ as test functions in the
relative entropy inequality (\ref{ec2}) and making use of a Gronwall type argument. This has been done in detail
in \cite[Section 6]{EF101} in the case of a viscous fluid satisfying the Navier-Stokes-Fourier system. However, the same arguments
can be used to handle the inviscid case provided the solutions are uniformly bounded on the existence interval.
\eR

\bRemark{ec3}
As the proof of Theorem \ref{Tec1} is based on the relative entropy inequality (\ref{ec2}), the conclusion remains valid if we replace
the internal energy \emph{equation} (\ref{i3}) by the entropy \emph{inequality} (\ref{i5b}) as long as we require (\ref{i5a}).
\eR

\subsection{Infinitely many dissipative solutions}

Apparently, the stipulation of the total energy balance (\ref{i5a}) eliminates the non-physical solutions obtained
in Theorem \ref{Tsec1}, at least in the case of \emph{regular} initial data. As we will see, the situation changes if we consider
non-smooth initial data, in particular the initial velocity field $\vu_0$ belonging only to $L^\infty(\Omega;R^3)$. Our final goal
is the following result.

\Cbox{Cgrey}{

\bTheorem{ec2}
Let $T> 0$ be given. Let the initial data $\vr_0$, $\vt_0$ be given, satisfying
\bFormula{ec4}
\vr_0 , \vt_0 \in C^2(\Omega),\ \vr_0(x) > \underline{\vr} > 0,\ \vt_0(x) > \underline{\vt} > 0 \ \mbox{for any}\ x \in \Omega.
\eF

Then there exists a velocity field $\vu_0$,
\[
\vu_0 \in L^\infty(\Omega;R^3),
\]
such that the problem (\ref{i1} - \ref{i5}) admits infinitely many dissipative (weak) solutions in $(0,T) \times \Omega$, with the initial
data $[\vr_0, \vt_0, \vu_0]$.
\eT

}

\bRemark{ec4}
As we shall see below, the solutions obtained in the proof of Theorem \ref{Tec2} enjoy the same regularity
as those in Theorem \ref{Tsec1}, in particular, the equation of continuity (\ref{i1}) as well as the internal energy balance
(\ref{i3}) are satisfied pointwise (a.a.) in $(0,T) \times \Omega$.
\eR

\bRemark{eeeec5}
In general, the initial velocity $\vu_0$ depends on the length of the existence interval $T$. 
See Section \ref{CCL} for more discussion concerning possible extension of the solutions to $[0, \infty)$.
\eR

The remaining part of this section is devoted to the proof of Theorem \ref{Tec2} that may be viewed as an extension of
the results of Chiodaroli \cite{Chiod} and De Lellis-Sz\' ekelyhidi \cite{DelSze3} to the case of a heat conducting fluid.

\subsubsection{Suitable initial data}

Following the strategy of \cite{DelSze3} our goal is to identify suitable initial data $\vu_0$ for which
the associated (weak) solutions of the momentum equation dissipate the kinetic energy. In contrast with \cite{DelSze3}, however,
we have to find the initial data for which the kinetic energy decays \emph{sufficiently} fast in order to compensate the associated production of heat.

The velocity field $\vc{v} = \vr \vu$ we look for will be \emph{solenoidal}, in particular, we focus on the initial data satisfying
\[
\Div (\vr_0 \vu_0) = 0.
\]
This assumption simplifies considerably the ansatz introduced in Section \ref{refor}, specifically,
\[
\vr = \tilde \vr = \vr_0(x), \ \vc{v} = \vr \vu ,\ \Div \vc{v} = 0, \Psi \equiv 0;
\]
whence the problem reduces to solving
\bFormula{ec5}
\partial_t \vc{v} + \Div \left( \frac{\vc{v} \otimes \vc{v}}{\tilde \vr} \right) + \Grad \left(  \tilde \vr \vt [ \vc{v}]
 - \frac{2}{3} \chi \right)  = 0,\ \Div \vc{v} = 0, \ \vc{v}(0, \cdot) = \vc{v}_0, \ \Div \vc{v}_0 = 0,
\eF
for a suitable spatially homogeneous function $\chi = \chi(t)$.

Mimicking the steps of Section \ref{reduc} we introduce the quantity
\bFormula{ec6}
\Ov{e}[\vc{v}] = \chi - \frac{3}{2} \tilde \vr \vt[\vc{v}].
\eF
As the anticipated solutions satisfy
\[
\frac{1}{2} \frac{ |\vc{v}|^2 }{\tilde \vr}  = \Ov{e}[\vc{v}],
\]
the energy of the system reads
\bFormula{AA}
E (t) = \intO{ \left( \frac{|\vc{v}|^2}{2 \tilde \vr} + \frac{3}{2} \tilde \vr \vt[\vc{v}] \right)(t, \cdot) } =
\chi(t) |\Omega| = \chi(t).
\eF
Consequently, in accordance with the construction procedure used in Section \ref{reduc}, it is enough to find a suitable
\emph{constant} $\chi$ and the initial velocity field $\vc{v}_0$ such that
\[
\Div \vc{v}_0 = 0,\
E_0 = \intO{ \left( \frac{|\vc{v_0}|^2}{2 \tilde \vr_0} + \frac{3}{2} \tilde \vr_0 \vt_0 \right) } = \chi,
\]
and the associated space of subsolutions $X_0$ defined in (\ref{i15}) (with $\Grad \Psi = 0$) is non-empty. This is the objective of the remaining part of this
section.

\subsubsection{Dissipative data for the Euler system}

Similarly to (\ref{i18}), we introduce the set of subsolutions
\bFormula{i18A}
X_{0, e}[T_1,T_2] = \Bigg\{ \vc{v} \ \Big| \ \vc{v} \in  C^1((T_1,T_2) \times \Omega;R^3)
\cap C_{\rm weak}([T_1,T_2]; L^2(\Omega;R^3))
\eF
\[
 \vc{v} \ \mbox{satisfies (\ref{i14}) with some}\ \tn{U} \in C^1((T_1,T_2) \times \Omega ; R^{3 \times 3}_{{\rm sym},0}),
\]
\[
\left.
\frac{3}{2} \lambda_{\rm max} \left[ \frac{ \vc{v}  \otimes \vc{v} }{\tilde \vr} - \tn{U} \right] < e
\ \mbox{for}\ t \in (T_1,T_2),\ x \in \Omega \right\},
\]
where $e \in C([T_1,T_2] \times \Omega)$.

The following result may be seen as an extension of \cite[Proposition 5]{DelSze3}:

\bLemma{A2} Suppose that $\vc{v} \equiv \vc{v}_0 (x)$, together with the associated field $\tn{U}_{\vc{v}} \equiv 0$,
belong to the set of subsolutions $X_{0,e}[0,T]$.

Then for any $\tau \in (0,T)$ and any $\ep > 0$, there exist $\Ov{\tau} \in (0,T)$, $|\tau - \Ov{\tau}| <\ep$ and $\vc{w} \in X_{0,e}[\Ov{\tau},T]$, such that
\bFormula{AAA1}
\frac{1}{2} \frac{|\vc{w}(\Ov{\tau}, \cdot)|^2}{\tilde \vr} = e(\Ov{\tau}, \cdot) ,
\eF
\[
\vc{w} \equiv \vc{v} , \ \tn{U}_{\vc{w}} \equiv 0 \ \mbox{in a (left) neighborhood of}\ T.
\]
\eL

\bRemark{ec5}
Note that, thanks to (\ref{AAA1}),
\[
\vc{w}(t, \cdot) \to \vc{w}(\Ov{\tau}, \cdot) \ \mbox{(strongly) in}\ L^2(\Omega;R^3)\ \mbox{as}\ t \to \Ov{\tau} +.
\]

\eR

\bRemark{ec6}
The result is probably not optimal; one should be able, with greater effort, to show the same conclusion with
$\Ov{\tau} = \tau$.
\eR

\bProof

We construct the function $\vc{w}$ as a limit of a sequence $\{ \vc{w}_k \}_{k=1}^\infty \subset X_{0,e}[0,T]$,
\[
\vc{w}_k \to \vc{w} \ \mbox{in}\ C_{\rm weak} ([0,T]; L^2(\Omega;R^3)),
\]
where $\vc{w}_k$ will be obtained recursively, with the starting point
\[
\vc{w}_0 = \vc{v} \equiv \vc{v}_0 ,\ \tau_0 = \tau,\ \ep_0 = \ep.
\]
More specifically, we construct the functions $\vc{w}_k$, together with
$\tau_k$, $\ep_k$, $k=1,\dots$ satisfying:

\begin{itemize}

\item

\bFormula{AA3-}
\vc{w}_{k} \in X_{0,e}[0,T],\ {\rm supp}[\vc{w}_k - \vc{w}_{k-1}] \subset (\tau_{k - 1} - \ep_k, \tau_{k - 1} + \ep_k),\
 \mbox{where}\ 0< \ep_k < \frac{\ep_{k-1}}{2};
\eF

\item

\bFormula{AA2-}
d(\vc{w}_k, \vc{w}_{k-1}) < \frac{1}{2^k}, \ \sup_{t \in (0,T)} \left| \intO{ \frac{1}{\tilde \vr} (\vc{w}_k - \vc{w}_{k-1}) \cdot \vc{w}_{m} } \right| < \frac{1}{2^k}
\ \mbox{for all}\ m=0,\dots, k-1,
\eF
recalling that $d$ is the metric induced by the weak topology of the Hilbert space $L^2(\Omega;R^3)$.

\item
there exists $\tau_k$,
\[
\tau_{k} \in (\tau_{k-1} - \ep_k,
\tau_{k-1} + \ep_k)
\]
such that
\bFormula{AA2}
\intO{ \frac{1}{2} \frac{|\vc{w}_{k}|^2}{\tilde \vr} (\tau_k, \cdot) }  \geq
\intO{ \frac{1}{2} \frac{|\vc{w}_{k-1}|^2}{\tilde \vr} (t, \cdot) }  + \frac{\lambda}{\ep^2_k} \alpha_k^2
\eF
\[
\geq \intO{ \frac{1}{2} \frac{|\vc{w}_{k-1}|^2}{\tilde \vr} (\tau_{k-1}, \cdot) }  + \frac{\lambda}{2 \ep^2_k} \alpha_k^2 \ \mbox{for all}\ t \in (\tau_{k-1} - \ep_k, \tau_{k - 1} + \ep_k),
\]
where
\[
\alpha_k = \int_{\tau_{k-1} - \ep_k}^{\tau_{k-1} + \ep_k} \intO{ \left( e - \frac{1}{2} \frac{ |\vc{w}_{k-1}|^2 } {\tilde \vr} \right) } \ \dt > 0,
\]
and $\lambda > 0$ is constant independent of $k$.
\end{itemize}

Supposing we have already constructed $\vc{w}_0, \dots, \vc{w}_{k - 1}$ we find $\vc{w}_{k}$ enjoying the properties (\ref{AA3-} - \ref{AA2}).
To this end, we first compute
\[
\alpha_k = \int_{\tau_{k-1} - \ep_k}^{\tau_{k-1} + \ep_k} \intO{ \left( e - \frac{1}{2} \frac{ |\vc{w}_{k-1}|^2 } {\tilde \vr} \right) } \ \dt
\ \mbox{for a certain} \ 0 < \ep_k < \frac{\ep_{k-1}}{2}
\]
and observe that
\[
\frac{\alpha_{k}}{2 \ep_{k}} = \frac{1}{2 \ep_{k}} \int_{\tau_{k-1} - \ep_{k}}^{\tau_{k-1} + \ep_{k}} \intO{ \left( e - \frac{1}{2} \frac{ |\vc{w}_{k - 1}|^2 } {\tilde \vr} \right) } \ \dt \to \intO{ \left( e - \frac{1}{2} \frac{ |\vc{w}_{k-1}|^2 } {\tilde \vr} \right)(\tau_{k-1}) } > 0
\ \mbox{for}\ \ep_{k} \to 0
\]
as $\vc{w}_{k-1}$ is smooth in $(0,T)$.

Consequently, by the same token, we can choose $\ep_k > 0$ so small that
\bFormula{odhad}
\frac{1}{2 \ep_k} \int_{\tau_{k-1} - \ep_{k} }^{ \tau_{k-1} + \ep_{k}}
\intO{ \frac{1}{2} \frac{|\vc{w}_{k-1}|^2}{\tilde \vr}  } \ \dt + \frac{\Lambda(\tilde \vr, \|e \|_{L^\infty(0,T) \times \Omega)}) }{4 \ep_{k}^2} \alpha_{k}^2
\eF
\[
\geq
\intO{ \frac{1}{2} \frac{|\vc{w}_{k-1}|^2}{\tilde \vr}(t,\cdot)  }  + \frac{\Lambda(\tilde \vr, \|e \|_{L^\infty(0,T) \times \Omega)})}{8 \ep_{k}^2} \alpha_{k}^2
\]
\[
\geq \intO{ \frac{1}{2} \frac{|\vc{w}_{k-1}|^2}{\tilde \vr}(\tau_{k-1},\cdot)  }  + \frac{\Lambda(\tilde \vr, \|e \|_{L^\infty(0,T) \times \Omega)})}{16 \ep_{k}^2} \alpha_{k}^2 \ \mbox{for all}\ t \in (\tau_{k-1} - \ep_k, \tau_{k-1} + \ep_k),
\]
where $\Lambda(\tilde \vr, \|e \|_{L^\infty(0,T) \times \Omega)}) > 0$ is the universal constant from Lemma \ref{LP2}.

Applying Lemma \ref{LP2} in the form specified in Remark \ref{Rw4} we obtain a function $\vc{w}_k \in X_{0,e}$
such that
\[
{\rm supp}[ \vc{w}_{k} - \vc{w}_{k-1} ] \subset (\tau_{k-1} - \ep_{k}, \tau_{k-1} + \ep_{k} ),
\]
\bFormula{aA1}
d(\vc{w}_{k}, \vc{w}_{k-1} ) < \frac{1}{2^{k}}, \ \sup_{t \in (0,T)} \left| \intO{ \frac{1}{\tilde \vr} (\vc{w}_{k} - \vc{w}_{k-1}) \cdot \vc{w}_{m} } \right| < \frac{1}{2^k},\ m=0,\dots, k-1
\eF
and
\bFormula{odhad2}
\int_{\tau_{k-1} - \ep_{k} }^{ \tau_{k-1} + \ep_{k}} \intO{ \frac{1}{2} \frac{|\vc{w}_{k}|^2}{\tilde \vr}  } \ \dt
\geq
\int_{\tau_{k-1} - \ep_{k} }^{ \tau_{k-1} + \ep_{k}}
\intO{ \frac{1}{2} \frac{|\vc{w}_{k-1}|^2}{\tilde \vr}  } \ \dt + \frac{\Lambda(\tilde \vr, \|e \|_{L^\infty(0,T) \times \Omega)})}{2 \ep_{k}} \alpha_{k}^2,
\eF
where we have applied Jensen's inequality to the last integral in (\ref{P1Da}).

Finally, the relations (\ref{odhad}), (\ref{odhad2}) yield (\ref{AA2}) with some $\tau_k \in (\tau_{k-1} - \ep_k, \tau_{k-1} + \ep_k)$,
$\lambda = \Lambda / 16$.

Now, by virtue of (\ref{AA2-}), there is $\vc{w}$ such that
\bFormula{conv}
\vc{w}_k \to \vc{w} \ \mbox{in} \ C_{\rm weak}([0,T]; L^2(\Omega; R^3)).
\eF
Moreover,  (\ref{AA3-}) implies {\bf (i)}
\[
\tau_k \to \Ov{\tau} \in (0,T) ,\ | \Ov{\tau} - \tau | < \ep;
\]
{\bf (ii)} for any $\delta > 0$ there is $k = k_0(\delta)$ such that
\bFormula{AA5}
\vc{w}(t, \cdot) = \vc{w}_{k} (t, \cdot) = \vc{w}_{k_0}(t, \cdot) \ \mbox{for all} \ t \in (0, \Ov{\tau} - \delta) \cup (\Ov{\tau} + \delta, T),
k \geq k_0.
\eF
In particular, (\ref{AA5}) yields
\[
\vc{w} \in X_{0,e}[\Ov{\tau},T], \ \mbox{and}\ \vc{w} \equiv \vc{v}, \ \tn{U}_{\vc{w}} \equiv 0 \ \mbox{in a (left) neighborhood of}\ T.
\]

Next, in view of (\ref{AA2}),
\bFormula{AA6-}
\intO{ \frac{1}{2} \frac{ |\vc{w}_{k-1} |^2 }{\tilde \vr} (t, \cdot) } \nearrow Y \ \mbox{uniformly for}\ t \in (\tau_{k-1} - \ep_k, \tau_{k-1} + \ep_k),
\eF
therefore
\bFormula{AA6}
\frac{\alpha_k}{\ep_k} = \frac{1}{\ep_k} \int_{\tau_{k-1} - \ep_k}^{\tau_{k-1} + \ep_k} \intO{ \left( e - \frac{1}{2} \frac{ |\vc{w}_{k-1}|^2 } {\tilde \vr} \right) } \ \dt \to 0;
\eF
whence, finally,
\bFormula{AA7}
\intO{ \frac{1}{2} \frac{ |\vc{w}_k|^2 }{\tilde \vr} (\Ov{\tau} , \cdot) } \nearrow \intO{ e(\Ov{\tau}, \cdot) }.
\eF

Combining (\ref{AA7}) with (\ref{AA2-}), (\ref{conv}) we get
\[
\vc{w}_k (\Ov{\tau}, \cdot) \to \vc{w}(\Ov{\tau}, \cdot) \ \mbox{in}\ L^2(\Omega; R^3)
\]
which implies (\ref{AAA1}). Indeed we have
\[
\intO{ \frac{1}{\tilde \vr} | \vc{w}_n - \vc{w}_m |^2 (\Ov{\tau}, \cdot) }
\]
\[
=
\intO{ \frac{1}{\tilde \vr} |\vc{w}_n |^2 (\Ov{\tau}, \cdot) } - \intO{ \frac{1}{\tilde \vr} |\vc{w}_m |^2 (\Ov{\tau}, \cdot) }
- 2 \intO{ \frac{1}{\tilde \vr} \left( \vc{w}_n - \vc{w}_m \right) \cdot \vc{w}_m (\Ov{\tau}, \cdot) }
\ \mbox{for all} \ n > m,
\]
where, by virtue of (\ref{AA2-}),
\[
\intO{ \frac{1}{\tilde \vr} \left( \vc{w}_n - \vc{w}_m \right) \cdot \vc{w}_m (\Ov{\tau}, \cdot) }
= \sum_{ k = 0}^{n-m-1} \intO{ \frac{1}{\tilde \vr} \left( \vc{w}_{k+1} - \vc{w}_k \right) \cdot \vc{w}_m (\Ov{\tau}, \cdot) }
\to 0 \ \mbox{for}\ m \to \infty.
\]
\qed

\subsubsection{Construction of suitable initial data for the Euler-Fourier system}

Fixing $\vr_0$, $\vt_0$ satisfying (\ref{ec4}) and $\vr = \tilde \vr \equiv \vr_0$ we can use
(\ref{i12}) to deduce that there is a constant $\Ov{\vt}$ depending only on $[\vr_0, \vt_0]$ such that
\bFormula{A2}
|\vt [\vc{v}]| \leq \Ov{\vt}  ,\ \mbox{whence}\
\frac{3}{2} p(\tilde \vr, \Ov{\vt}[\vc{v}]) < \Ov{P} \ \mbox{on the whole interval}\ [0,T],
\eF
with $\Ov{P}$
independent of $\vc{v}$.

Next, we estimate the difference
$\vt - \vt_0$ satisfying the equation
\[
\tilde \vr \partial_t (\vt - \vt_0) + \vc{v} \cdot \Grad (\vt - \vt_0) - \frac{2}{3} \Delta (\vt - \vt_0) = - \vc{v} \cdot \Grad \vt_0 + \frac{2}{3} \Delta \vt_0 + \frac{2}{3}
\vt \vc{v} \cdot \frac{\Grad \tilde \vr}{\tilde \vr}.
\]
Consequently, using (\ref{A2}) and the comparison principle, we deduce that
\bFormula{A3}
|\vt [\vc{v}](t, \cdot) - \vt_0 | \leq c\left( 1 + \|\vc{v} \|_{L^\infty((0,T) \times \Omega;R^3)} \right) t \ \mbox{for all}\ t \in [0,T].
\eF

We take $\vc{v}_0 \in C^1(\Omega)$, $\Div \vc{v}_0 = 0$, and a constant $\chi_0$ in such a way that
\bFormula{A4}
\frac{3}{2} \lambda_{\rm max} \left( \frac{\vc{v}_0 \otimes \vc{v}_0 }{\tilde \vr} \right) < \chi_0 - \frac{3}{2} \vr_0 \vt_0.
\eF
Moreover, for any $\Ov{\chi} > 2 \chi_0$, $K > 0$ given, it is easy to construct a function $\chi \in C[0,T]$ such that
\begin{itemize}
\item
\[
\chi(0) = \chi(T) = \chi_0, \ \chi(t) > \chi_0 \ \mbox{for all}\ t \in (0,T), \ \max_{t \in (0,T)} \chi(t) = \Ov{\chi};
\]
\item
there is $\tau \in (0,T)$ and $\ep > 0$ such that
\bFormula{form}
\chi(\Ov{\tau}) - \chi_0 > \frac{\Ov{\chi}}{2}  , \ \chi(t) < \chi(\Ov{\tau}) - K (t - \Ov{\tau}) \ \mbox{for all} \ t \in (\Ov{\tau} , T)
\ \mbox{whenever}\ |\tau - \Ov{\tau} | < \ep.
\eF
\end{itemize}

Consequently, we have
\[
\vc{v} \equiv \vc{v}_0 \in X_{0,e}[0,T] \ \mbox{(with}\ \tn{U} \equiv 0 )
\]
provided
\[
e(t,x) = \chi(t) - \frac{3}{2} \vr_0 \vt_0.
\]

Applying Lemma \ref{LA2}, we find a function $\vc{w} \in X_{0,e}[\Ov{\tau}, T]$, with the corresponding field $\tn{U}_{\vc{w}}$, such that
\[
\frac{1}{2} \frac{|\vc{w}(\Ov{\tau}, \cdot)|^2}{\tilde \vr} = \chi(\Ov{\tau}) - \frac{3}{2} \vr_0 \vt_0 > \frac{\Ov{\chi}}{2} + \chi_0 - \frac{3}{2} \vr_0 \vt_0,
\]
\[
\vc{w} \equiv \vc{v}_0, \ \tn{U}_{\vc{w}} = 0 \ \mbox{in a (left) neighborhood of}\ T,
\]
and
\[
\frac{1}{2} \frac{| \vc{w} |^2}{\tilde \vr} < \frac{3}{2} \lambda_{\rm \max} \left( \frac{\vc{w} \otimes \vc{w}}{\tilde \vr} - \tn{U}_{\vc{w}} \right)
< \chi(t) - \frac{3}{2} \vr_0 \vt_0
\]
\[
\leq \chi(\Ov{\tau}) - \frac{3}{2} \vr_0 \vt_0 - K(t - \Ov{\tau}), \ t \in (\Ov{\tau}, T].
\]

Denoting $\vc{w}_0 = \vc{w}(\Ov{\tau}, \cdot)$ and
shifting everything to the origin $t = 0$, we infer that there is a function $\vc{w} \in X_{0,e}(0,T)$, with
the following properties:
\begin{itemize}
\item
\bFormula{A5}
\vc{w}(0, \cdot) = \vc{w}_0, \ \frac{1}{2} \frac{|\vc{w}_0 |^2}{\tilde \vr} = \chi(\Ov{\tau}) - \frac{3}{2} \vr_0 \vt_0,
\ \vc{w}(T, \cdot) = \vc{v}_0,
\eF
\item
\bFormula{A6}
e(t,x) = \left\{ \begin{array}{c} \chi(\Ov{\tau})  - \frac{3}{2} \vr_0 \vt_0 - K t ,\ t \in \left[0, \frac{ \chi(\Ov{\tau}) - \chi_0}{K} \right]
\\ \\ \chi_0 - \frac{3}{2} \vr_0 \vt_0 \ \mbox{for}\ t \in \left[\frac{ \chi(\Ov{\tau}) - \chi_0}{K} , T \right] . \end{array} \right.
\eF

\end{itemize}

Similarly to (\ref{i15}), we introduce the set $X_0$, together with the function
\[
\Ov{e}[\vc{v}] = \chi(\Ov{\tau})  - \frac{3}{2} \vr_0 \vt[\vc{v}].
\]
Our ultimate goal is to show that the function $\vc{w}$, introduced in (\ref{A5}), belongs to $X_0$ as long as we conveniently fix the parameters
$\Ov{\chi}$, $K$. To this end, it is enough to show that $e$, defined through (\ref{A6}), satisfies
\bFormula{A7}
e(t, x) < \Ov{e}[\vc{w}] =  \chi(\Ov{\tau}) - \frac{3}{2} \vr_0 \vt[\vc{w}] \ \mbox{for all}\ t \in (0,T].
\eF

For $t \in \left( 0,\frac{ \chi(\Ov{\tau}) - \chi_0}{K} \right]$, this amounts to showing
\[
\frac{3}{2} \left( \vr_0 \vt[\vc{w}] - \vr_0 \vt_0 \right) < Kt, \  t \in \left( 0,\frac{ \chi(\Ov{\tau}) - \chi_0}{K} \right],
\]
which follows from (\ref{A3}) provided $K = K(\Ov{\chi})$ is taken large enough.

Next, for $t \in \left[\frac{ \chi(\Ov{\tau}) - \chi_0}{K} ,T \right]$, we have to check that
\[
\frac{3}{2} \left( \vr_0 \vt[\vc{w}] - \vr_0 \vt_0 \right) < \frac{\chi(\Ov{\tau})}{2} ,\ t \in \left[\frac{ \chi(\Ov{\tau}) - \chi_0}{K} ,T \right],
\]
which follows from (\ref{A2}), (\ref{form}) provided we fix $\Ov{\chi} = \Ov{\chi}(\vr_0, \vt_0)$ large enough.

Having found a suitable subsolution we can finish the proof
of Theorem \ref{Tec2} exactly as in Section \ref{reduc}.

\section{Concluding remarks}
\label{CCL}

\begin{itemize}

\item
In this paper, we focused exclusively on the physically relevant $3D-$ setting. The reader will have noticed that exactly the same results
may be obtained also in the $2D-$case. Note, however, that the method does not apply to the $1D-$system as the conclusion of
Lemma \ref{LP1} is no longer available.

\item
Theorem \ref{Tec2} obviously applies to the larger class of dissipative solutions for which the internal energy balance is replaced by
the entropy \emph{inequality}
\bFormula{ccl1}
\partial_t \left( \vr \log \left( \frac{\vt^{3/2}}{\vr} \right) \right) + \Div \left(\vr \log \left( \frac{\vt^{3/2}}{\vr} \right) \vu \right) - \Div \left( \frac{ \Grad \vt }{\vt} \right) \geq  \frac{|\Grad \vt|^2}{\vt^2}.
\eF
Moreover, we could even construct dissipative solutions with an ``artificial'' entropy production satisfying
(\ref{ccl1}) with \emph{strict} inequality and, at the same time, conserving the total energy.
On the other hand, a criterion based on maximality of the entropy production could be possibly used to identify a class of physically relevant solutions.

\item The conclusion of Theorem \ref{Tec1} can be extended to the time interval $[0, \infty)$ by means of continuation. Indeed
we can take the function $h$ in (\ref{sec3aa}) such that $h(T) = 0$; whence
\[
\tilde \vr(T, \cdot) = \vr_0.
\]
Moreover,
as pointed out in Remark
\ref{Rrem1}
\[
\vt(t, \cdot) \in [L^p(\Omega), W^{2,p}(\Omega) ]_{\alpha} \ \mbox{for all}\ t \in [0,T],
\]
therefore we can apply Theorem \ref{Tec1} recursively on the time intervals $[nT, (n+1)T]$, n=1,\dots

A similar extension of Theorem \ref{Tec2} seems possible but technically more complicated.

\end{itemize}

\centerline{
\textbf{Acknowledgment}}

\medskip

The authors would like to thank Camillo De Lellis for fruitful discussions about the problem.


\begin{thebibliography}{10}

\bibitem{AL1}
T.~Alazard.
\newblock Low {M}ach number flows and combustion.
\newblock {\em SIAM J. Math. Anal.}, 38(4):1186--1213 (electronic), 2006.

\bibitem{AL}
T.~Alazard.
\newblock Low {M}ach number limit of the full {N}avier-{S}tokes equations.
\newblock {\em Arch. Rational Mech. Anal.}, {\bf 180}:1--73, 2006.

\bibitem{Amann1}
H.~Amann.
\newblock Nonhomogeneous linear and quasilinear elliptic and parabolic boundary
  value problems.
\newblock In {\em Function spaces, differential operators and nonlinear
  analysis ({F}riedrichroda, 1992)}, volume~{\bf 133} of {\em Teubner-Texte
  Math.}, pages 9--126. Teubner, Stuttgart, 1993.

\bibitem{BiaBre}
S.~Bianchini and A.~Bressan.
\newblock Vanishing viscosity solutions of nonlinear hyperbolic systems.
\newblock {\em Ann. of Math. (2)}, {\bf 161}(1):223--342, 2005.

\bibitem{BRESSAN}
A.~Bressan.
\newblock {\em Hyperbolic systems of conservation laws. The one dimensional
  Cauchy problem.}
\newblock Oxford University Press, Oxford, 2000.

\bibitem{Chiod}
E.~Chiodaroli.
\newblock A counterexample to well-posedness of entropy solutions to the
  compressible {E}uler system.
\newblock 2012.
\newblock Preprint.

\bibitem{D4}
C.~M. Dafermos.
\newblock {\em Hyperbolic conservation laws in continuum physics}.
\newblock Springer-Verlag, Berlin, 2000.

\bibitem{DelSze3}
C.~De~Lellis and L.~Sz{\'e}kelyhidi, Jr.
\newblock On admissibility criteria for weak solutions of the {E}uler
  equations.
\newblock {\em Arch. Ration. Mech. Anal.}, {\bf 195}(1):225--260, 2010.

\bibitem{DelSze}
C.~De~Lellis and L.~Sz{\'e}kelyhidi, Jr.
\newblock The {$h$}-principle and the equations of fluid dynamics.
\newblock {\em Bull. Amer. Math. Soc. (N.S.)}, {\bf 49}(3):347--375, 2012.

\bibitem{EF101}
E.~Feireisl.
\newblock Relative entropies in thermodynamics of complete fluid systems.
\newblock {\em Discr. and Cont. Dyn. Syst. Ser. A}, {\bf 32}:3059--3080, 2012.

\bibitem{FeNo6}
E.~Feireisl and A.~Novotn{\' y}.
\newblock {\em Singular limits in thermodynamics of viscous fluids}.
\newblock Birkh{\" a}user-Verlag, Basel, 2009.

\bibitem{FeiNov10}
E.~Feireisl and A.~Novotn{\' y}.
\newblock Weak-strong uniqueness property for the full
  {N}avier-{S}tokes-{F}ourier system.
\newblock {\em Arch. Rational Mech. Anal.}, {\bf 204}:683--706, 2012.

\bibitem{Krylov}
N.~V. Krylov.
\newblock Parabolic equations with {VMO} coefficients in {S}obolev spaces with
  mixed norms.
\newblock {\em J. Funct. Anal.}, {\bf 250}(2):521--558, 2007.

\bibitem{LiuTP}
T.~P. Liu.
\newblock Admissible solutions of hyperbolic conservation laws.
\newblock {\em Mem. Amer. Math. Soc.}, 30(240):iv+78, 1981.

\bibitem{MulSve}
S.~M{\"u}ller and V.~{\v{S}}ver{\'a}k.
\newblock Convex integration for {L}ipschitz mappings and counterexamples to
  regularity.
\newblock {\em Ann. of Math. (2)}, {\bf 157}(3):715--742, 2003.

\bibitem{Serr3}
D.~Serre.
\newblock Local existence for viscous system of conservation laws: {$H^s$}-data
  with {$s>1+d/2$}.
\newblock In {\em Nonlinear partial differential equations and hyperbolic wave
  phenomena}, volume~{\bf 526} of {\em Contemp. Math.}, pages 339--358. Amer.
  Math. Soc., Providence, RI, 2010.

\bibitem{Serr4}
D.~Serre.
\newblock The structure of dissipative viscous system of conservation laws.
\newblock {\em Phys. D}, {\bf 239}(15):1381--1386, 2010.

\bibitem{Shn}
A.~Shnirelman.
\newblock Weak solutions of incompressible {E}uler equations.
\newblock In {\em Handbook of mathematical fluid dynamics, {V}ol. {II}}, pages
  87--116. North-Holland, Amsterdam, 2003.

\bibitem{WIL}
C.~H. Wilcox.
\newblock {\em Sound propagation in stratified fluids}.
\newblock Appl. Math. Ser. 50, Springer-Verlag, Berlin, 1984.

\end{thebibliography}

\end{document}